\documentclass[11pt, a4]{amsart}

\usepackage[latin1]{inputenc}
\usepackage{amsfonts}
\usepackage{latexsym}
\usepackage{amssymb}
\newtheorem{thm}{Theorem}[section]
\newtheorem{prop}[thm]{Proposition}
\newtheorem{lem}[thm]{Lemma}
\newtheorem{cor}[thm]{Corollary}
\newtheorem{defi}[thm]{Definition}
\newtheorem{rmq}[thm]{Remark}

\def\preuve{\smallskip\goodbreak{\it Proof.~~--~\kern.3em}
\ignorespaces}
\def\qedbox{$\square$}
\def\qed{\ifmmode\qedbox\else\unskip\ \hglue0mm\hfill
\qedbox\smallskip\goodbreak\fi}

\newcommand{\C}{\ensuremath{\mathbb C}}
\newcommand{\R}{\ensuremath{\mathbb R}}

\newcommand{\p}{\ensuremath{\mathbb P}}
\newcommand{\diff}{\ensuremath{\mbox{\textit{Diff}\/}}} 
\newcommand{\gal}{\ensuremath{\mbox{\textit{Gal}\/}}} 
\newcommand{\aut}{\ensuremath{\mbox{\textit{Aut}\/}}}  
\newcommand{\hol}{\ensuremath{\mbox{\textit{Hol}\/}}}   
\newcommand{\inv}{\ensuremath{\mbox{\textit{Inv}\/}}}   

\newcommand{\F}{\ensuremath{\mathcal F}}

\begin{document}

%%%%%%%%%%%%
%%  TITRE
%%%%%%%%%%%

\title[The Galoisian envelope...]{The Galoisian envelope of a germ of foliation:\\ the quasi-homogeneous case}

\author{\sc E. PAUL }

\address{
E. Paul\\
Laboratoire Emile Picard UMR 5580 \\
Institut de Mathématiques de Toulouse\\
Universit\'{e} Paul Sabatier \\
118 route de Narbonne \\
31062 Toulouse cedex 9, France.
}

\email{paul@math.ups-tlse.fr}

\date{\today}
\subjclass{}
\keywords{}
\thanks{Work partially supported by the EU GIFT project (NEST- Adventure Project
n° 5006)}

\maketitle

\textit{Abstract.} {\footnotesize We give geometric and algorithmic criterions in order to have of a proper Galois closure for a codimension one germ of quasi-homogeneous foliation. We recall this notion recently introduced by B. Malgrange, and describe the Galois envelope of a group of germs of analytic diffeomorphisms. The geometric criterions are obtained from transverse analytic invariants, whereas the algorithmic ones make use of formal normal forms.}

%\tableofcontents

\section*{Introduction}
There are several notions of integrability for a system of differential equations. Most of them are related to the existence of a sufficient number of first integrals for the solutions of the system. These definitions differ each other on the additional properties required for this family of invariants functions. We can separate them into two types:\\
- conditions between the first integrals: one may ask commutativity conditions for the Poisson bracket, or relax such a condition;\\
- conditions on the nature of these functions: rational, meromorphic or multivalued functions in some ''reasonable'' class of transcendence.\\
The main methods for proving non integrability (analytical methods, Ziglin method or Morales-Ramis method) are based on the linearization of the system around a particular solution. Therefore they only deliver sufficient criterions on non integrability, using for the last mentioned method \textit{linear} differential Galois theory. 

In order to investigate the second type of condition, and --in the future-- to get necessary and sufficient conditions for integrability, we have to consider the system in the whole, which suggests to consider a  \textit{non linear} differential Galois theory. The first attempts in this direction was done by J. Drach and E. Vessiot. More recently, B. Malgrange introduced in \cite{MAL1} (see also the introductive version \cite{MAL2}) a ''Galois closure'' for any dynamical system, namely the smallest D-groupoid which contains the solutions of the system. Roughly speaking, a D-groupoid is a system of partial differential equations whose local solutions satisfy groupoid conditions outside an analytic codimension one set. They are not strict Lie groupoid, in order  to deal with singular systems. As a matter of introduction to this notion, we shall describe in the first section the Galois closure of a group of germs of analytic diffeomorphisms at the origin of \C . 

Each D-groupoid admits a D-algebra obtained by the linearization of its equations along the identity solutions. The local solutions of this linear differential system are stable under the Lie bracket outside of a codimension one analytic set. The Galois envelope of a singular analytic foliation \F\ is the smalllest D-groupoid \gal(\F) whose D-algebra contains the germs of tangent vector fields to \F . It is a proper one if it doesn't coincide with the whole groupoid \aut(\F) obtained by writing the equations of invariance of the foliation under a local diffeomorphism. In this case  --which is not the general case--, its solutions satisfy an additional differential relation, and we shall say that the foliation is Galois reducible. 

For a local codimension one singular foliation defined by a holomorphic one-form $\omega$, this reducibility property is equivalent to the existence of a Godbillon-Vey sequence of finite length for $\omega$ (at most three): there exists a finite sequence of lenght at most three of meromorphic one forms $\omega_0$, $\omega_1$, and $\omega_2$ such that $\omega_0$ is an equation of the foliation and
$$d\omega_0=\omega_0\wedge\omega_1,\ d\omega_1=\omega_0\wedge\omega_2,\ d\omega_2=\omega_1\wedge\omega_2.$$
This fact was described in manuscripted notes by B. Malgrange \cite{MAL0}, and then has been extensively proved by G. Casale in \cite{CAS3} with some different arguments. In particular, the transverse rank of \gal\F\ (i.e. the order of its transverse local expression) is also the minimal lenght of a Godbillon-Vey sequence for $\F$. Finally, G. Casale proved in \cite{CAS0} that this Godbillon-Vey condition is also equivalent to the existence of first integrals for the foliation with a particular type of transcendence which belongs to a Darboux or Liouville or Riccati type differential extension, according to the transverse rank of the Galois closure. These different points of view on the Galois reducibility admit a generalization for higher codimension foliations: see \cite{CAS4} for Painlevé 1 foliation.

\bigskip

In the present paper we shall only deal with codimension one foliations. Therefore, we expect the existence of at most one first integral, and we only have to discuss the second type of integrability condition: the existence of such a first integral in a given class of transcendency. The previous discussion allows us to reformulate the integrability problem as following:  \textit{give necessary and sufficient criterions for the Galois reducibility of a germ of codimension one foliation.}
We present an answer to this problem in the following context: \F\ is defined by a vector field $X=X_h+\cdots$ where the ''initial'' hamiltonian vector field $$X_h=\frac{\partial h}{\partial y}\frac{\partial}{\partial x}-\frac{\partial h}{\partial x}\frac{\partial}{\partial y}$$ is quasi-homogeneous with to respect to $R=p_1x\frac{\partial}{\partial x}+p_2y\frac{\partial}{\partial y}$ ($p_1$, $p_2$ positive integers): $R(h)=\delta h$, $\delta=\deg_R(h)$. The dots means terms of higher quasihomogeneous degree. We furthermore require that $h$ has an isolated singularity (with Milnor number $\mu$) and that $X$ still keep invariant the analytic set $h=0$. Therefore, $X$ is a logarithmic vector field for the polar set $h=0$, and we have:
$$X=aX_h+bR,\ a\in \mathcal{O}_2, b\in \mathcal{O}_2,\ a(0)=1$$
with $deg_R(bR)>\deg_R(X_h).$
The restriction to this class of foliation is motivated by the two following reasons:\\
- the desingularization of these foliations by blowing up's is ''simple'': it is similar to the one of the quasi-homogeneous function $h$: the exceptional divisor is only a chain of projective lines and all the pull-back of the irreducible components of $h$ --excepted the axis if they appear in $h$-- meet the same ''principal'' projective line $C$. \\
- in this class of foliations, we have at our disposal \textit{formal normal forms} which give us complete formal invariants: see \cite{P2}. 

This will allow us to give two different types of criterions for the Galois reducibility of $\F$ : a geometric one which is related to the holonomy of the principal component $C$ of the desingularized foliation, and an algorithmic one which directly holds on the normalized formal equation of the foliation. For the first one, let us denote $\hol(\F)$ the holonomy group of the principal component $C$ for the desingularized foliation. This is an analytic invariant of \F\ (in fact, this ''transverse invariant'' is also a complete invariant in this quasi-homogeneous context: see \cite{YG}). We prove in theorem (\ref{critere.geom}) the following result : 

\medskip
\textbf{Theorem 1. }\textit{The Galois groupoid of the germ of quasi-homogeneous foliation \F\ is a proper one if and only if the Galois envelope of $\hol({\F})$ is a proper one.}
\medskip

This theorem reduces the initial problem to the determination of the Galois closure of a  subgroup $G$ of \diff(\C,0), which is described in the first section (theorem \ref{groupes.galoisiens}). The main argument in the proof of this theorem is an extension of the equation which define the Galois closure of $\hol(\F)$ to the whole exceptional divisor. This is possible, since the elements of the holonomy group of $C$ are solutions of this equation and therefore keep it invariant. This proof suggests that even in non quasi-homogeneous cases, these criterions for the Galois reducibility will only depend on the transverse structure of the foliation.

Theorem 1 is not an explicit criterion since in general, we can't compute the invariant $\hol(\F)$. In order to get an algorithmic criterion, we recall in section 3 the formal normal forms for this class of foliations. Notice that in general these models are divergent models. The radial component of these normal forms make appear a collection $\mathcal{L(F)}$ of $\mu$ formal one-variable vector fields, and it turns out that this collection (up to a commun conjugacy) is a complete invariant for the formal class of $\F$. It must be surprising to try to charaterize the Galois reducibility of \F\ using only formal invariants. Nevertheless, we can perform it according to the two following facts:\\
- if a foliation is Galois reducible, then its formal normal form is a convergent one;\\
- if the foliation \F\ is a ''non exceptional'' one (see \cite{CM}), then there exits a convergent conjugacy between \F\ and its model. \\
Clearly, for exceptional foliations, we need an additional condition on the analytic class of \F , which is not yet an algorithmic one. The central result of this work is the following theorem which summarize theorem \ref{Red.gal.for}, corollary \ref{Red.gal.nonex} and theorem \ref{Red.gal.ex}:

\medskip
\textbf{Theorem 2.} \textit{If the quasi-homogeneous foliation \F\ is a non exceptional one, the Galois closure of \F\ is proper if and only if the explicit invariant $\mathcal{L}(\F)$ generates a finite dimensional Lie algebra. In this case, this one is always of dimension one, and the foliation is at most Liouvillian.\\
\indent  If the quasi-homogeneous foliation \F\ is an exceptional one, the Galois closure of \F\ is proper if and only if the explicit invariant $\mathcal{L}(\F)$ is a finite dimensional Lie algebra, and the analytic invariants of \F\ are of ''unitary'' or ''binary'' type. In this case, the foliation will be a Liouvillian one (for unitary invariants), or of Riccati type, (for binary invariants).}

\medskip

We shall recall in the first section the definition of unitary or binary invariants which is a terminology introduced by J. Ecalle. The first part of the theorem is an extension of a result of F. Loray and R. Meziani for nilpotent singularities \cite{LM}, while the second one is an extension of a theorem of G. Casale for reduced singularities \cite{CAS3}. Notice that in the local context, the Galois reducible foliations which are not Liouvillian are very rare.

\bigskip

Clearly, the relationship between the algorithmic invariant $\mathcal{L}(\F)$ and the geometric one $\hol({\F})$ has a transcendental nature since the first one is directly obtained from the differential equation whereas the second one is related to the solutions of this equation. Nevertheless, for Galois reducible foliations we can describe this relationship: it reduces to the exponential map of the one-variable vector fields of $\mathcal{L}(\F)$. In order to check this fact it is more convenient to consider an equivalent data to $\hol({\F})$: the relative holonomy of \F\ with respect to its initial part defined by $X_h$ (see section 4). 
 
Finally, we conclude this paper with a list of open questions related to the present results.

\bigskip

\section{The Galois closure of a subgroup of \diff(\C,0)}

\medskip

Let $\Delta$ be a disc around 0 in \C. We first recall the list of all the D-groupoids on $\Delta$ (see \cite{MAL0} and \cite{CAS1}). We denote $(x,y,y_1,y_2,\cdots y_k)$ the coordinates for the space of $k$-jets of maps from $\Delta$ to itself.

\begin{thm}\label{liste1} The differential ideal of a $D$-groupoid on $\Delta$ is generated by a meromorphic equation of one of the five types:
\begin{enumerate}
\item $D$-groupoids of order zero: they are generated by an equation of the form: $h(x)-h(y)=0$ where $h$ is a holomorphic function on $\Delta$. We denote them: $G_0(h)$. 
\item $D$-groupoids of order one: they are generated by an equation of the form: $\eta(y)(y_1)^n-\eta(x)=0$ where $n$ is an integer, and
  $\eta$ a meromorphic function on $\Delta$. We denote them $G_1^n(\eta)$.
\item $D$-groupoids of order two: they are generated by an equation of the form: $\mu(y)y_1+\frac{y_2}{y_1}-\mu(x)=0$
where $\mu$ is meromorphic on $\Delta$. We denote them $G_2(\mu)$.
\item $D$-groupoids of order three: they are generated by an equation of the form: 
  $\nu(y){y_1}^2+2\frac{y_3}{y_1}-3\left(\frac{y_2}{y_1}\right)^2-\nu(x)=0$
  where $\nu$ is meromorphic on $\Delta$. We denote them $G_3(\nu)$.
\item The $D$-groupoid of infinite order $G_\infty$ defined by the trivial equation $0=0$, whose solutions are the whole sheaf $\aut(\Delta)$.
\end{enumerate}
\end{thm}

The Galois envelope of a subgroup $G$ of \diff(C,0) is the smallest D-groupoid in the previous list which admits all the elements $g$ of $G$ as solutions. Clearly, the existence of a proper Galois closure of finite order $k$, only depends on the analytic class of $G$.
The Galois envelope for a monogeneous subgroup generated by $g$ is the Galois envelope of $g$ itself, since all the iterates of $g$ will also satisfy the same equation, by composition or inversion stability. The Galois closure $\gal(g)$ of $g$ is given by the two following results, see B. Malgrange  \cite{MAL0}, and G. Casale (\cite{CAS1}). Let $\alpha=g'(0)$. If $\alpha$ is an irrational number, then $g$ is formally linearizable. We have:

\begin{prop}\label{lin}
A formally linearizable diffeomorphism has a proper Galois closure if and only if it is  an analytically linearizable diffeomorphism. In this case, its Galois closure is a rank one D-groupoid.
\end{prop}

If $\alpha$ is a rational number, $g$ is a resonant diffeomorphism, and there exists an integer $q$ such that $g^q$ is tangent to the identity. The following lemma
\begin{lem}\label{puissance} \cite{CAS1} For all non vanishing  integer $q$, $\gal(g)=\gal(g^q)$.
\end{lem}
\noindent reduces the study to the case $\alpha=1$. Any diffeomorphism tangent to the identity to an order $k$ is conjugated via a formal series to a normal form $g_N$ which is the exponential of the vector field 
$\frac{x^{k+1}}{1+\lambda x^k}\frac{d}{dx}$. Following the description of J. Martinet and J.P. Ramis, we obtain a complete analytic invariant $\inv(g)$ of $g$ by the following construction (see \cite{MR1}). Using 2$k$ sectorial normalizations, one can prove that the space of the orbits of $g$ is obtained by gluing $2k$ bipunctured Riemann spheres $(S_i,0,\infty)$ 
with local diffeomorphisms $\varphi_i^0: (S_{i-1},0)\rightarrow (S_i,0)$ and $\varphi_{i}^\infty: (S_i,\infty)\rightarrow (S_{i+1},\infty)$. The collection \inv($g$) of these gluing maps up to global automorphisms on each $(S_i,0,\infty)$ is a complete analytic invariant of $g$. This invariant \inv(g) is \textit{unitary}  if there exists a positive integer $p$ such that the gluing maps $\varphi_i^0$ are of the form $z\mapsto z/(1+a_iz^p)^{1/p}$  and  $\varphi_i^\infty$: $u\mapsto u$ at infinity ($u=1/z$). It is a \textit{binary} one when the gluing maps are alternatively of the form
$z\mapsto z/(1+a_iz^p)^{1/p}$ in 0 and $u\mapsto u/(1+b_iu^p)^{1/p}$ at infinity. We have:

\begin{thm}(see \cite{CAS1})\label{Galois(g)}
Let $g$ be an element of \diff(\C,0) tangent to the identity. The Galois closure $ \gal(g)$
is proper of rank two (resp. three) if and only if its analytic invariant $\inv(g)$ is a unitary one (resp. a binary one).
\end{thm}

\begin{rmq}\label{theta(g)}
The proof of this theorem make use of the following result.
Let $\theta_g$ be the formal vector field such that $g=\exp\theta_g$ (there is existence and unicity of such vector field, and its multiplicity at 0 is greater or equal to 2). The diffeomorphism $g$ is a solution of a D-groupoid if and only the vector field $\theta_g$ is a formal solution of its D-Lie algebra.
\end{rmq}

We now discuss the Galois closure of a subgroup $G$ of \diff(\C,0) generated by $g_1,\cdots g_\mu$. Let $\Theta$ (resp. $\widehat{\Theta}$)
be the Lie algebra of one variable holomorphic (resp. formal) vector fields which vanish at the origin: $\theta=(a_kz^{k}+\cdots)d/dz$. Recall that (see for example \cite{LOR}):

\begin{lem}\label{dim-finie} A subalgebra $\mathcal G$ of $\Theta$ (resp. $\widehat{\Theta}$) is a finite dimensional one if and only if $\mathcal G$ is at most of dimension two. Furthermore, such a Lie algebra is always a solvable one, and if the multiplicity $k$ of each element of $\mathcal G$ is greater or equal to two, then its dimension is at most one.
\end{lem}

Notice that such a result, and thus the following proposition, cannot be generalized in a global situation, in which there exist three dimensional Lie algebras of one variable vector fields which are not solvable ones.

\begin{prop}\label{gal.res} If the subgroup $G$ of \diff(\C,0) has a proper Galois closure, then $G$ is a solvable group.\end{prop}

\begin{preuve}
Let $G_1$ be the subgroup of $G$ whose elements are tangent to the identity map. If $G_1$ is trivial, then $G$ is abelian since the first derivative group $[G,G]$ of $G$ is contained in $G_1$, and we are done. Therefore, we suppose that $G_1$ is non trivial. For each element $g$ of $G_1$, let $\theta_g$ be the element of $\widehat{\Theta}$ such that $g=\exp\theta_g$. From remark \ref{theta(g)}, the Lie algebra $\mathcal{L}(G_1)$ generated by these vector fields is included in the solutions of the D-Lie algebra of the Galois closure of $G$, and is a finite dimensional one. Therefore, from the previous lemma, its dimension is one, and there exists a vector field $\theta$ and constants $c_g$ such that for all $g$ in $G_1$, $g=\exp c_g\theta$. This proves that $G_1$ is an abelian group. Since $[G,G]$ is contained in $G_1$, the group $G$ is a solvable one.
\end{preuve}

The converse of Proposition (\ref{gal.res}) is false: all the monogeneous subgroups are abelian, but from (\ref{Galois(g)}), outside of the unitary or binary cases, they don't have a proper Galois closure. The Galois reducibility is not only an algebraic property of the group $G$.

\bigskip

We shall recall the formal classification of the solvable subgroups of \diff(\C,0) (see \cite{P2} or \cite{CM}). We denote $G_1$ the subgroup of $G$ of its elements tangent to the identity map. We have:

\medskip
-   the group $G$ is formally linearizable if and only if $G_1$ is the trivial group;

\medskip
- every solvable non linearizable group $G$ is formally conjugated to a group $G_N$ of the following type:
$$G_N=\{g_{\lambda,t}=\lambda\exp t\frac{z^{k+1}}{1+\alpha z^k}\frac{d}{dz},\
  \lambda\in\Lambda, \ t\in T\}$$
where $\Lambda$ is a multiplicative subgroup of $\C^*$ and $T$ is an additive subgroup of $\C$. Furthermore, $G_N$ is abelian if and only if
$\Lambda$ is a group of $k$-roots of 1. If $G_N$ is not an abelian group, the residue $\alpha$ vanishes, and the elements of $G_N$ are obtained by lifting homographies fixing 0 with the ramification $z\mapsto z^k$.

\medskip
-  Following the terminology of D. Cerveau and R. Moussu \cite{CM}, $G$ is an \textit{exceptional} subgroup of \diff(\C,O) if
$G_1$ is  monogeneous. In particular, they are solvable groups. These authors prove that, among the non linearizable groups, the non exceptional groups are exactly the rigid ones: the formal classification coincide with the analytic one. One should say that an exceptional group is a unitary or binary one when $G_1$ is generated by a unitary or binary element.
\begin{thm}\label{groupes.galoisiens} The only subgroups of \diff(\C,0) which have a proper Galois closure are:

(1) the analytically linearizable groups;

(2) the non exceptional solvable groups;

(3) the exceptional unitary groups;

(4) the exceptional binary groups.

Furthermore, the rank of their D-envelope is at most one in case (1), at most two in cases (2) and (3), and at most three in case (4). 
\end{thm}

We call \textit{Liouvillian} group every subgroup of \diff(\C,0) whose Galois closure is at most of rank two, and \textit{Riccatitian} group every subgroup of \diff(\C,0) whose Galois closure is at most of rank three. In the present local situation, the Riccatitian non Liouvillian groups are very rare: their class is restricted to the (non empty!) set defined by $(4)\setminus (3)$.

\bigskip

\begin{preuve}
We first check that these groups have a proper Galois closure:

\medskip
(1) Since the existence of a proper Galois closure of finite order $k$ only depends on the analytic class of $G$, it suffices to consider a group of linear diffeomorphisms. They keep invariant the differential form $dx/x$ and therefore satisfy the differential equation $xy_1-y=0$ which is, according to the notations of (\ref{liste1}), the equation of the $D$-groupoid $G_1^1(1/x)$. Remark that this is only an upper bound of $\gal(G)$: for example, if $G$ is a group of periodic rotations, they keep invariant an holomorphic function $h$ and $\gal(G)=G_0(h)$.

\medskip
(2) The formal model $G_N$ of a solvable group is Liouvillian. Indeed, the differential form $\omega=(1+\alpha
   x^k)/x^{k+1}dx$ is invariant by each element $f_{\lambda , t}$ of $G_N$ up to a multiplicative constant
   $c_{\lambda, t}$. Therefore, each element of $G_N$ satisfy
$a(y)y_1=c_{\lambda , t}a(x)$, where $a$ is the coefficient of $\omega$. Derivating these equations, each element of $G$ is a solution of the same equation
$$a(x)a(y)y_2+a'(y)a(x)y_1^2-a'(x)a(y)y_1=0$$
where $a'$ is the derivative of $a$ with respect to $x$. This is the equation of the rank two D-groupoid $G_2(a'/a)$.
The same previous remark holds: this is only an upper bound of the Galois closure of $G_N$: if $G_N$ is abelian, its elements all satisfy the rank one equation $a(y)y_1-a(x)=0$ of $G_1^1(a)$. Now, if $G$ is a non exceptional group, by rigidity, it is analytically conjugated to $G_N$, and still have a proper Galoisian envelope of rank at most two.

\medskip
(3) and (4): Let $G$ be an exceptional group and let $g_1$ be a generator of the monogeneous group $G_1$, which is supposed to be unitary or binary. From (\ref{Galois(g)}), $G_1$ has a proper closure of rank two or three with equation $E=0$. If $G$ is not equal to $G_1$, we know from proposition 2 of \cite{CM} that $G$ is generated by $g_1$ and a second resonant element $g_2$. If $g_1$ is tangent to the identity at order $k$, the normal form of $G$ described by \cite{CM} shows that $g_2^{2k}$ belongs to $G_1$ and therefore $g_2^{2k}=g_1^l$ for some integer $l$. With lemma (\ref{puissance}), we conclude that $g_2$ also belongs to the Galois closure of $G_1$, and finally, $\gal(G)=\gal(G_1)$.

\medskip

On the converse, we now suppose that $G$ has a proper Galois envelope. If $G_1$ is a trivial group, then $G$ is formally linearizable and from proposition (\ref{lin}) we conclude that $G$ is of type (1). If $G_1$ is non trivial, we know from (\ref{gal.res}) that $G$ is a  solvable group. Either it is a non exceptional one, and $G$ is of type (2), or it is an exceptional one: $G_1$ is generated by an element $g_1$. Since the Galois closure of this one is non trivial, we know from theorem (\ref{Galois(g)}) that $g_1$ and thus $G$ is of type (3) or (4).
\end{preuve}

\bigskip

\section{A geometric criterion for Galois reducibility}

\medskip

We first recall general facts on the Galois reducibility for singular holomorphic foliations. Let $\mathcal F$ be a singular holomorphic foliation of codimension $k$ on a $n$-dimensional holomorphic manifold of $M$. Following the definition of B. Malgrange \cite{MAL1}, the Galois groupoid of $\mathcal F$ is its D-envelope, i.e. the smallest D-groupoid $\gal({\mathcal F})$ wich is ''admissible'' for the foliation: its D-Lie algebra contains the tangent vector fields. The Galois groupoid of $\mathcal F$ is always contained in the D-groupoid $\aut({\mathcal F})$ of the germs of diffeomorphisms which keep invariant the foliation. We shall say that $\mathcal F$ is Galois reducible if its Galois envelope is proper: $\gal({\mathcal F})\neq \aut({\mathcal F})$. This property only depends on the analytic class of the foliation, and is invariant by blowing up or blowing down transformations.

If $U$ is an open set in $M$ on which the foliation is trivializable by tangent-transverse coordinates $(s,t)$, $s=(s_1,\cdots s_{n-k})$, $t=(t_1,\cdots t_k)$, the local ideal of $\gal({\mathcal F})$ can be generated by equations (see \cite{CAS3}):
\begin{eqnarray}\label{st}\frac{\partial T_j}{\partial z_i}=0,\ E_i(t,T,\cdots  \frac{\partial^{\mid\alpha\mid} T}{\partial t^{\alpha}})\end{eqnarray}
where $E_i$ are the equations of a D-groupoid on the $k$-dimensional polydisc $t(U)$. The rank of this local transverse groupoid doesn't depend on the local chart \cite{CAS3}: this is the \textit{transverse rank} of $\gal({\mathcal F})$.

We now suppose that $\mathcal F$ is a codimension one foliation on a polydisc $\Delta$ in $(\C^n,0)$, defined by a one-form $\omega$ which satisfies the Frobenius condition. We may suppose that the singular locus is at least a codimension two analytic set. From (\ref{liste1}), the transverse rank of $\mathcal F$ can only get the values 0, 1, 2, 3 or $\infty$, the finite values corresponding to the proper cases. A Godbillon-Vey sequence for $\omega$ is a sequence of meromorphic one-forms $\omega_n$ such that
$$d\omega=\omega\wedge\omega_1,\ \ d\omega_1=\omega\wedge\omega_2,\cdots$$
$$d\omega_i=\omega\wedge\omega_{i+1}+\sum_{j=1}^{i}\binom{i}j\omega_j\wedge\omega_{i-j+1}$$
A Godbillon-Vey sequence of lenght $l>1$ is a Godbillon-Vey sequence such that $\omega_i=0$, $i\geq l$. A Godbillon-Vey sequence of lenght 1, is a Godbillon-Vey sequence of lenght 2, such that $\omega_1=p^{-1}df/f$ for an integer $p$: $f^{1/p}$ is an integrating factor of $\omega$. The existence of a Godbillon-Vey sequence of lenght $l$ only depends on the foliation defined by $\omega$. We have (see \cite{MAL0} and \cite{CAS3}):

\begin{thm}\label{GV}
The foliation $\mathcal F$ has a Godbillon-Vey sequence of lenght $l$ with $l\leq 3$ if and only if the transverse rank of its Galois groupoid is at most $l$.
\end{thm}

Furthermore, G. Casale has proved in \cite{CAS0} that the existence of a proper Galois envelope for $\mathcal F$ is also equivalent to the existence of a transcendental first integral which belongs to a particular type of extension, namely a meromorphic, Darboux, Liouvillian or Riccatician type, according to the values $l=0,1,2$ or 3 of the transverse rank of $\gal({\mathcal F})$. Therefore, in each case, one should call the foliation with the same terminology.

\bigskip

If $L$ is a leaf of $\mathcal F$, and if $\hol(L)$ is the image of its holonomy representation, then all its elements are solutions of the local ideal of $\gal({\mathcal F})$. Indeed, for any loop $\gamma$ which represents an element of $\pi_1(L,m)$, we can cover $\gamma$ by trivializing open sets $U_1, \cdots U_p$ such that the transverse coordinate on $U_i$ is an analytic extension of the previous one. With this choice, the change of local coordinates are tangent to the foliation and therefore are solutions of $\gal({\mathcal F})$. By the stability under composition, the change of coordinates between $U_p$ and $U_1$ is a solution of $\gal({\mathcal F})$. In particular, its transverse component --which is the holonomy representation of $\gamma$-- is a solution of the local expression of $\gal({\mathcal F})$. From this remark, and since the existence of a proper Galois envelope is an invariant property under birational maps, we obtain 

\begin{prop}\label{gal-hol} If $\mathcal F$ has a proper Galois envelope, then any holonomy group of $\mathcal F$ or of any foliation $\widetilde{\mathcal F}$ obtained from $\mathcal F$ by blowing up's has a proper Galois closure whose rank is at most the transverse rank of $\gal({\mathcal F})$.
\end{prop} 

We shall prove that for the present class of quasi-homogeneous germs of foliations, we have a converse of this statement. In order to do this, we consider the desingularization process of \F : see \cite{SEI} or \cite{MM}. For a quasi-homogeneous foliation which is a perturbation of the foliation defined by $h=0$, extending an argument of \cite{CM}, one can prove that the desingularization process is the same as the one of $dh$, namely: the exceptional divisor is a chain of projective lines which are invariant for the desingularized foliation; all the strict pull back of each component of $h=0$ different from the axis are transverse to the same projective line $C$: we call it the principal one. One can check that $C$ is also the space of the values for the meromorphic first integral $x^{p_2}/y^{p_1}$ of the quasi-radial vector field $R$. The singularities on $C$ are the different values corresponding to each branch of $X$, and 0, $\infty$, which are the intersections with other components. If $x$ or $y$ occurs in the decomposition of $h$, their pullback by the composition of blowing up's  is a line transverse to the end components of the chain. All the reduced singularities are resonnant saddles (not necesseraly linearizable), since their linear part is obtained by the local expression of the desingularization of $dh/h$. The projective holonomy of $\mathcal F$ is the holonomy of the principal component $C$ of the desingularized foliation $\widetilde{\mathcal F}$. We denote \hol(\F) the image of this representation: this is a subgroup of $\diff(\C,0)$ defined up to a conjugacy (the choice of a transverse on which we realize the holonomy group).
The following result is announced in \cite{YG}, and proved for cuspidal singularities in \cite{M}:

\begin{thm}\label{class-anal} Two quasi-homogeneous germs of foliations $ \F_1$ and $\F_2$ are analytically equivalent if and only if $\hol(\F_1)$ is conjugated to $\hol(\F_2)$.
\end{thm}

\noindent The easier following result can be proved independently: 

\begin{thm}\label{critere.geom} The Galois envelope of the germ of quasi-homogeneous foliation \F\ is a proper one if and only if the Galois envelope of \hol(\F) is a proper one.
\end{thm}

\begin{preuve}
If the Galois envelope of \F\ is a proper one, the same holds for \hol(\F) from proposition (\ref{gal-hol}). We now suppose that \hol(\F) has a Galois envelope of finite rank given by an equation $E=0$ of type (0), (1), (2) or (3) in the list given by (\ref{liste1}). Let $(s_0,t_0)$ be a local system of tangent-transverse coordinates on an open set $U_0$ around a regular point $m$ in the principal component of $\widetilde{\mathcal F}$, and let $T$ be the transversal $s_0=s_0(m)$. As above, we can extend $E=0$ to a local equation $E_0$ of a D-groupoid on $U_0$ whose transverse expression is $E=0$ and is admissible for the foliation setting: 
$$\frac{\partial T}{\partial s}=0,\ E(t,T,\cdots  \frac{\partial^{k} T}{\partial t^{k}})=0.$$
We can extend this D-groupoid along a path $\gamma$ by covering this path with open sets $U_\alpha$, $\alpha=0,\cdots n$ with local systems $(s_\alpha, t_\alpha)$: the first equation is preserved by a foliated change of coordinates, and the second one $E_\alpha=0$ is extended on $U_\beta$ by $\psi^*_{\alpha\beta}E_\alpha$ where $t_\beta=\psi_{\alpha\beta}(t_\alpha)$. If $\gamma$ is a loop, this analytic extension coincide at the end of $\gamma$ with the initial groupoid: indeed, the composition of the transition maps $\psi_{\alpha\beta}$ is the holonomy map of $\gamma$ and we know that this one is a solution of the Galois envelope, and therefore keep invariant $E_\alpha=0$. By this way, we get an extension of the D-groupoid $E_0=0$ to the smooth part of the principal component $C$. 
Now, we can extend this groupoid to a neibourhood of each reduced singularity on $C$, from a result of Guy Casale: see proposition (5.2) in \cite{CAS3}. 
Let $C'$ be an adjacent component to $C$ and $p$ a regular point near from $C\cap C'$. One can choose local generators of the groupoid in $p$ which are still writen under the previous adapted form. Furthermore, the local holonomy of $C'$ around $C\cap C'$ is a solution of this groupoid. From the previous description of the exceptional divisor, $C'$ gets at most two singularities, and the fundamental group of the complement of its singularities is generated by one element. Therefore we can extend the groupoid along $C'$ and inductively to the whole divisor.
\end{preuve}

Such a type of argument can be used to prove that if $\hol(\F_1)$ is conjugated to $\hol(\F_2)$, then this conjugacy gives a local conjucagy around $m$ for the desingularized corresponding foliations, whose \textit{transverse} expression can be extended to the whole divisor. The main difficulty in theorem (\ref{class-anal}) is to prove that for quasihomogeneous  foliations there is no \textit{tangent} obstruction to construct a global conjugacy along the divisor. Here, the existence of a proper Galois envelope -or of a Godbillon-Vey sequence: see \cite{P2} for the Liouvillian case- only involves transverse obstructions, and thus are easier to obtain.

\bigskip

\section{An algorithmic criterion for Galois reducibility}

\medskip

We want to test the Galois reducibility by making use of formal normal forms for the germs of quasi-homogeneous foliations $\F_X$ defined by the vector fields:
$$X=aX_h+bR,\ a\in \mathcal{O}_2, b\in \mathcal{O}_2,\ a(0)=1, \deg(bR)>\deg(X_h).$$
In the general situation, both normal forms and conjugacies are formal objects. This will only give a criterion of \textit{formal} Galois reducibility. We can consider two definitions for the \textit{formal} Galois reducibility of an analytic foliation:\\
(i-) There exists an analytic foliation $\F'$ which is formally conjugated to \F\ and Galois reducible;\\
(ii-) The foliation $\F$ admits a formal finite Godbillon-Vey sequence.\\
Clearly the first one implies the second one, by taking with the formal conjugacy the pull back of the Godbillon-Vey sequence of $\F'$ given by (\ref{GV}). We first choose the second definition here, since we deal with formal models. But finally, it turns out that, for our class of foliations, both definitions coincide (see remark (\ref{fnc}) below): when this criterion of formal Galois reducibility holds, we shall obtain convergent final normal forms. Therefore, if we are in a non exceptional (or ''rigid'') case, the conjugacy will also converge, and we shall obtain an algorithmic criterion for analytic Galois reducibility.

\bigskip
 We first recall the construction of these normal forms and introduce the related complete formal invariant, obtained in \cite{P2}. They generalize the normal forms of the cuspidal case ($h=y^2-x^3$) described in \cite{ZOL} and \cite{LOR1}. We split it into two steps:

\medskip
\textit{First step: prenormalization.} It is based on the following general lemma.
Let $M$ be a submodule of the $\widehat{\mathcal O}_n$-module of formal vector fields at the origin of $\C^n$, endowed with a graduation, and stable under the Lie bracket (in the present case, $M$ is the module of logarithmic vector fields, endowed with the quasi-homogeneous degree induced by $R$). Let $X=X_0+\cdots $ be a perturbation of the initial quasi-homogeneous vector field $X_0$ of degree $\delta_0$ by higher order terms.

\begin{lem}\label{prenorm}\cite{P2} Let $B$ be the image of the operator $[X_0,\cdot]$ in $M$, and $W$ a complement space of $A=B+\widehat{\mathcal O}_nX_0$ in $M$. There exits a vector field $Y$ in $W$, a formal diffeomorphism $\Phi $ and a formal unity $u$ such that
$\Phi ^*X=u(X_0+Y).$
\end{lem}

Notice that if we want to classify the vector fields instead of the foliations (i.e. if we don't work up to a unity) the same statement holds with a complement of $B$ instead of $A$. This lemma reduces the first step to an appropriate choice of a submodule $W$ isomorphic to the quotient space $M/A$. Denote by $\mathcal I$ the ring of first integrals of the initial vector field $X_0$. The rich cases for normal forms occur when $\mathcal I$ doesn't reduce to the constants. In our case, $\mathcal{I}=\C[[h]].$ Clearly, the quotient $M/A$ is a $\mathcal I$-module, and one should naturally require the same property in our choice for $W$. In our present situation ($X_0=X_h$ with a quasi-homogeneous function $h$ which has an isolated singularity), we can compute the quotient $M/A$ (see \cite{P2} for details): this a free $\mathcal I$-module generated by the $\mu$ classes of vector fields $a_kR$, where $a_1,\cdots a_\mu$ is a monomial basis of the \C-vector space ${\mathcal O}_2 /J(h)$, and $J(h)$ is the jacobian ideal of $h$. This allows us to choose $W=X_h\oplus_{k=1}^{\mu}\mathcal I a_k\ R$, and from lemma (\ref{prenorm}), we have

\begin{thm}\label{etape1}
Let $X=aX _h +b R $ be a perturbation of $X _h$. There exits an element $(d_1,\cdots d_{\mu })$ of $\C[[h]]^{\mu}$, a formal diffeomorphism $\Phi $ which conjugate the foliation $\F_X$ to the foliation defined by the vector field
$$Y=X _h +\sum_{k=1}^{\mu }d_k(h){a_kR}.$$ 
Furthermore, we can require that this conjucagy is ''fibered'' with respect to $R$, i.e. is formally the exponential of a vector field proportional to $R$. Such a conjugacy keep invariant each trajectory of $R$.
\end{thm}

\textit{Second step: final reduction.} In the previous step, for a fixed complement space, there is no unicity of the prenormal form $Y$.  One can prove that the set of prenormal forms for $\F_X$ is the orbit of one of them under the action of a final reduction group of transformations of the following type: $\Phi=\exp b\cdot R$, with a formal coefficient $b$ in $\mathcal I$: $b=b(h)$. Such transformations satisfy the relation $h\circ \Phi=\varphi\circ h$ for a one variable formal diffeomorphism $\varphi$. In order to study the action of this final reduction group on the prenormal forms, it is convenient to introduce a modified expression of them. We shall make use of the two following remarks:

i- Setting $\alpha=h^{-\delta_0/\delta}$, he have $[\alpha X_h,R]=0$. The introduction of this multivalued coefficient will allow us to work with an abelian basis of logarithmic vector fields.

ii- Setting $r_i=\frac{\deg(\alpha a_i)}{\delta}$ we have 
$R(\alpha a_ih^{-r_i})=0$. This will allow us to work with coefficients which are constants for $R$.

Multiplying $Y$ with $\alpha$, and grouping coefficients in order to transform coefficients $a_i$ in constants $f_i$ for $R$ we obtain the following ''adapted'' prenormal forms:
\begin{eqnarray}\label{FNFmodifiee}
X_\alpha+ \sum_{i=1}^{\mu}f_i \delta_i(h)R
\end{eqnarray}
with $X_\alpha=\alpha X_h$, $f_i=\alpha a_ih^{-r_i}$ and $\delta_i=d_i(h)h^{r_i}$. 
By these two tricks, any element $\Phi$ of the final reduction group keep invariant
$X_\alpha$ and the coefficients $f_i$. Therefore we have
$$\Phi^*(X_\alpha+ \sum_{i=1}^{\mu}f_i \delta_i R)=X_\alpha+ \sum_{i=1}^{\mu}f_i \Phi^*(\delta_i(h) R).$$
The action of $\Phi$ over $\delta_i(h) R$ is given by
$$\Phi^*(\delta_i(h) R)= d_i\circ \varphi(h)\frac{\varphi(h))^{r_i+1}}{\varphi'(h)}\frac{R}{h}$$
where $\varphi$ is defined by $h\circ \Phi=\varphi\circ h$. This is the lifted action by $h$ of the action of $\varphi$ on the one-variable vector fields
$$\theta_i(z)=d_i(z)z^{r_i+1}\frac{d}{dz}.$$
Since $r_i=p_i/\delta$ for a positive integer $p_i$, we can uniformize these vector fields setting $t=z^{1/\delta}$ in
$$\theta_i(t)=\delta^{-1}d_i(t^\delta)t^{p_i+1}\frac{d}{dt}.$$
We may choose $\varphi$ -and therefore $\Phi$- in such a way that one of the vector fields $\theta_i$ is normalized under its usual normal form
$$\delta^{-1}\frac{t^{q_i+1}}{1+\lambda t^{q_i}}\frac{d}{dt}, \ \mbox{with }
q_i=\delta k_i+p_i$$
where $k_i$ is the multiplicity of each series $d_i$. Going back to the non adapted prenormal forms, we obtain the following final normal forms:

\begin{thm}\label{etape2}
Let $Y=X _h +\sum_{k=1}^{\mu }d_k(h){a_kR}$ be a prenormal form of $X$ and
$i$ an indice arbitrary chosen among $1,\cdots \mu$. There exits a diffeomorphism in the final reduction group which conjugate $Y$ to a normal form $Y_N$ in which the coefficient of indice $i$
is a rational function of $h$ of the following type: 
$$d^N_{i}(h)=\frac{h^m}{1+\lambda h^{m+n}}$$
where $\lambda$ is a complex number, and $m$, $n$ are positive integers.
\end{thm}

In fact, the previous argument gives rise to the following explicit formal invariant:

\begin{prop}\label{classification.formelle}
The family of the $\mu$ formal vector fields $\theta_i(t)$ up to a common conjugacy is a complete formal invariant for the foliation defined by $X$. We denote it $\mathcal{L}(\F).$
\end{prop}

Notice that as soon as $\mu$ is greater than two, we can't normalize simultaneously all the coefficients $d_i$ under a rational form. The final normal form is still a formal object. A result of M. Canalis and R. Schafke in the cuspidal situation ($h=y^2-x^3$) suggests that these final normal forms are defined by $k$-summable series in $t$: see \cite{CDS}. Nevertheless, the generalization of this fact, and the computation of the order $k$ is still an opened question. Furthermore, even if they are of the same nature (conjugacy class of $\mu$ one variable objects) the relationship between this algorithmic invariant $\mathcal{L}(\F)$ and the geometric one $\hol(\F)$ is not clear (it is of transcendental nature), excepted in the Galois reductive situations, in which we shall be able to specify it in the next section.

\bigskip
 
We now give a criterion of formal Galois reducibility, for the class of quasi-homogeneous foliations described in the introduction. 

\begin{thm}\label{Red.gal.for} The following propositions are equivalent:\\
(1) The foliation $\mathcal F$ is formally Galois reducible;\\
(2) The Lie algebra generated by the elements of $\mathcal{L}(\mathcal{F})$ is a finite dimensional one;\\
(3) The Lie algebra generated by the elements of $\mathcal{L}(\mathcal{F})$ is one dimensional;\\
(4) $\mathcal F$ is a formally Liouvillian foliation.
\end{thm}

\begin{preuve}
The equivalence between propositions (2) and (3) comes from Lemma (\ref{dim-finie}), since one can check that the multiplicity of each vector field $\theta_i$ is greater than one: this is a consequence of $\deg(bR)>\deg(X_h).$

We now prove the implication $(3) \Rightarrow (4)$. Let
$\theta$ and $c_k$ be a vector field and $\mu$ constants such that $\theta_k=c_k\theta$. The adapted normal form obtained in (\ref{FNFmodifiee}) is here of the following type:
$$X_\alpha+(\sum_{i=1}^{\mu}c_if_i)\delta(h) R.$$
The final reduction step normalize $\theta$ -and therefore here all the $\theta_i$- under its usual rational normal form. We obtain a \textit{convergent} normal form $X_N$ in the formal class of $\F_X$. In order to prove that the foliation $\mathcal F_N$ defined by $X_N$ is Galois reducible of order two, by theorem (\ref{GV}) we have to prove that there exits two logarithmic one-forms $\omega_N$ and $\omega_1$ such that $\omega_N$ define the foliation $\F_N$ and $\omega_1$ is a closed form such that $d\omega_N=\omega_N\wedge \omega_1.$ We consider the logarithmic one forms (for details on this dual point of view, see \cite{P2}):

$$\omega_h=\delta^{-1}\frac{dh}{h}=\frac{dx\wedge dy}{h}(\delta^{-1}X_h,\cdot), \ \ 
\omega_R=\frac{p_2ydx-p_1xdy}{h}=\frac{dx\wedge dy}{h}(\cdot,R).$$
Since $dx\wedge dy/h(\delta^{-1}X_h,R)=\delta^{-1}R(h)/h=1$, $\{\omega_R,\omega_h\}$ is a dual basis of $\{X_h,R\}$ for the pairing $(\omega,X)=\omega(X)$. Therefore, for any function $f$, we have $df=R(f)\omega_h+X_h(f)\omega_R$, and the one-form $a\omega_h-b\omega_R$ define the same foliation as $X=aX_h+bR$. Notice that $\omega_R$ is not a closed form, but dividing it with $\alpha=h^{-\delta_0/\delta}$, we have $d(\alpha^{-1}\omega_R)=0$: this is similar to the trick (i-). The foliation $\F_N$ is defined by 
$$\omega_N=\omega_h-\sum_{i=1}^{\mu}a_id_i(h)\omega_R=\omega_h-\sum_{i=1}^{\mu}f_i\delta_i(h)\frac{\omega_R}{\alpha}=
\omega_h-f_c\delta(h) \frac{\omega_R}{\alpha}$$
where $f_c=\sum_{i=1}^{\mu}c_if_i$ only depends on $c=(c_1,\cdots c_\mu)$. Since $R(f)=0$, we have $d(\omega_N/\delta(h))=0$, and the logarithmic derivative $\omega_1$ of $\delta(h)$ is a closed form which satisfies the Godbillon-Vey relation.

We now prove the main implication $(4) \Rightarrow (3)$. We shall give another proof of it in the next section.
If \F\ is formally Liouvillian then $\F_N$ have a (formal) Godbillon-Vey sequence of lenght two given by $\omega_N$, $\omega_1$, and  it suffices to prove that $\mathcal{L}(\F_N)$ is one dimensional. We can check that $\omega_1$ also keep invariant $X:\ h=0$, with simple poles along $X$ (for this last point, which is only formally true, see \cite{P1}). Therefore, $\omega_1$ is a closed logarithmic form and there exists two formal coefficients $\lambda$ and $\mu$ such that $\omega_1=\lambda\omega_h+\mu\omega_R$. We may suppose that $\omega_1=\lambda _0\omega_h$, where $\lambda_0$ in the residue of $\omega_1$ along $h=0$
, even if it means replacing $h$ with $h\circ\Phi$ and replacing the logarithmic basis with its pull back by $\Phi$: 
indeed, $\omega_1-\lambda_0\omega_h$ is a closed logarithmic form with vanishing residue and therefore, there exists a formal coefficient $g$ such that $\omega_1-\lambda_0\omega_h=\lambda_0dg$. Setting $u=\exp g$ we obtain: $\omega_1=\lambda_0d(uh)/(uh)$. Since $h$ is a quasihomogeneous function there exists a change of variable $\Phi$ such that $h\circ \Phi=uh$. Conjugating the Godbillon-Vey relation by $\Phi$ we normalize $\omega_1$ under the previous form, and $\Phi^*\omega_N$ is still normalized relatively to the new logarithmic basis.

\noindent Using the relations $df=R(f)\omega_h+X_h(f)\omega_R$ and  $d(\omega_R/\alpha)=0$, we obtain
$$d\omega_N=d(\omega_h-\sum_{i=1}^{\mu}f_i\delta_i(h)\frac{\omega_R}{\alpha})=- \sum_{i=1}^{\mu}f_iR(\delta_i(h))\omega_h\wedge\frac{\omega_R}{\alpha}.$$
Therefore, the Godbillon-Vey relation $ d\omega_N= \omega_N\wedge\omega_1$ is equivalent to
$$\lambda_0\ \sum_{i=1}^{\mu}f_i\delta_i(h)=-\sum_{i=1}^{\mu}f_i\delta'_i(h)\delta h.$$
where $\delta'$ is the derivative of this one-variable function. The decomposition of any element $b$ under a sum $\sum_{i=1}^{\mu}f_i\delta_i(h))$ or equivalently under a sum
$ \sum_{i=1}^{\mu}a_id_i(h))$ is unique. Indeed, the space of prenormal forms $W$ is isomorphic to the $\mathcal I$-module coker$(X_h)$, and this one is a free module over the basis given by the classes of $a_1,\cdots a_\mu$ (see \cite{P2}). Therefore the Godbillon-Vey equation is equivalent to the $\mu$ linear differential equations 
$$\delta h\delta'_i(h)=\lambda_0\delta_i(h),\ \ i=1\cdots \mu.$$
Since the functions $\delta_i(h)$ are solutions of the same one-dimensional first order linear differential equation, we have $\delta_i(h)=c_i\delta(h)$ for all $i$ in $\{1,\cdots \mu\}$.

Finally, we have to prove the non trivial implication of $(1) \Leftrightarrow (4)$, i.e.: any formally Galois reducible foliation is a formally Liouvillian one. This is essentially a consequence of (\ref{gal.res}). Indeed, if $\F$ has a proper Galois closure, we know that the same holds for $\hol(\F)$. From theorem (\ref{groupes.galoisiens}),  $\hol(\F)$
is a solvable group. According to \cite{P1}, this allows us to construct a formal Godbillon-Vey sequence of lenght two for the foliation. We summarize this construction: from Theorem (1.7) of \cite{P1}, a solvable subgroup of \diff (\C,0) admits a formal symmetry  i.e. a formal one variable vector field which is invariant up to a multiplicative constant by each element of the group. Evaluating $\omega$ on this symmetry, we obtain a local integrating factor whose logarithmic derivative $\omega_1$ satisfy the Godbillon-Vey relation. We can extend $\omega_1$ on the regular part of the principal component the exceptional divisor, since it is invariant by the holonomy of this component. Then, we extend it along the whole exceptional divisor with similar arguments as in the proof of (\ref{critere.geom}). 
\end{preuve}

\begin{rmq}\label{fnc} If the foliation \F\ is formally Galois reducible, then its final normal form is a convergent one. Indeed, 
if the Lie algebra generated by the elements of $\mathcal{L}(\mathcal{F})$ is one dimensional, then the action of the final reduction group will simultaneously normalize each coefficient $d_i(h)$ under a rational form. Therefore the final normal form of \F\ has a convergent expression. 
\end{rmq}

Following and extending the definition of D. Cerveau and R. Moussu in \cite{CM}, a quasi-homogeneous foliation \F\ is a non exceptional foliation if and only if its invariant $\hol(\F)$ is a non exceptional group. Two non exceptional holomorphic foliations which are formally conjugated are analytically equivalents: indeed, by \cite{P2}, we know that we can construct a conjugacy which is fibered with respect to $R$, and which is only a transversally formal one. Therefore, the restriction of such a transformation to any fiber of $R$ will define a conjugacy between the holonomy groups. Since they are non exceptional this conjugacy is a convergent one. 

\begin{cor}\label{Red.gal.nonex} We suppose that the foliation $\mathcal F$ is a non exceptional one. The following propositions are equivalent:\\
(1) The foliation $\mathcal F$ is Galois reducible;\\
(2) The Lie algebra generated by the elements of $\mathcal{L}(\mathcal{F})$ is a finite dimensional one;\\
(3) The Lie algebra generated by the elements of $\mathcal{L}(\mathcal{F})$ is one dimensional;\\
(4) $\mathcal F$ is a Liouvillian foliation. 
\end{cor}

\begin{preuve} The first implication $(1)\Rightarrow(2)$ (or (3)) comes from the corresponding implication in Theorem (\ref{Red.gal.for}). Since  $(4)\Rightarrow(1)$ is trivial, 
we only have to prove $(3)\Rightarrow(4)$. Let \F\ be a holomorphic foliation such that $\mathcal{L}(\mathcal{F})$ is one dimensional. Following the previous remark (\ref{fnc}), its final normal form is a convergent one and define a holomorphic foliation $\F_N$, which is Liouvillian. Since \F\ is a non exceptional foliation, the conjucacy beteween \F\ and $\F_N$ is a convergent one and \F\ is also a Liouvillian foliation. 
\end{preuve}

Notice that for non exceptional germs of foliations, there doesn't exist Riccatitian foliations which are not Liouvillian. Clearly, for an exceptional foliation, we need an additional criterion on the analytic class itself inside the formal one (they are all formally Liouvillian). This one is given by theorem (\ref{critere.geom}), and by the classification of the groups of diffeomorphisms with proper envelope (\ref{groupes.galoisiens}), and therefore is not yet an algorithmic one:

\begin{thm}\label{Red.gal.ex} 
An exceptional foliation \F\ has a proper Galois envelope if and only if the group \hol(\F) is an exceptional unitary or binary one.
\end{thm}

\bigskip

\section{Relationship between geometric and algorithmic invariants\\ for Galois reducible foliations}

\medskip

We introduced in section 2 the notion of projective holonomy, namely the holonomy of the principal component $C$ in the desingularization of the foliation. For explicit computations, the following  notion of ''relative holonomy'' is more efficient.
Let $m$ be a regular point of the desingularized  foliation on $C$ and $T$ the pull back of the fiber of $R$ corresponding to this value $m$. Any element of $\pi_1(C,m)$ can be lifted into a path from a point in $T$ in a leaf of the initial foliation $\F_h$ defined by $dh=0$. We consider the normal subgroup $\pi_1'(C,m)$ of $\pi_1(C,m)$ corresponding to the elements which can be lifted in \textit{loops} in the initial foliation: this is the kernel of the representation of the projective holonomy of $\F_h $. The \textit{relative representation of holonomy of} \F\ is the restriction of the projective holonmy
to $\pi_1'(C,m)$. We denote by $\hol(\F,\F_h)$ its image. This is a subgroup of the group $\diff_1(\C,0)$ of germs of diffeomorphisms tangent to the identity. We have:

%\bigskip

i-  The Galois closure of $\hol(\F,\F_h)$ is identical to the one of $\hol(\F)$;

ii- The class of conjugacy of $\hol(\F,\F_h)$ is still an analytic complete invariant for \F .

These two facts only hold for foliations \F\ which are a perturbation of $\F_h$. Indeed, in this case any element of $\hol(\F)$ is a resonnant one, and the statement (i-) is a consequence of proposition (\ref{puissance}). The second one can be deduced from theorem (\ref{class-anal}) by proving that the relative holonomy groups are conjugated if and only if the projective holonomy groups are conjugated. We don't give the details since we shall not make use of this result.

\bigskip

 The main interest of this holonomy is a more efficient presentation of $\pi_1'(C,m)$ interpreting its elements as \textit{horizontal classes of evanescent cycles}. Let us develop this point of view.
 We first remark that the quasi-radial vector field $R$ is a basic vector field for the initial foliation $\F_h$: from $R(h)=\delta h$, we deduce that its flow $\exp [\tau]R$ send the fiber   
$F_{z_0} :\ h=z_0$ on the fiber $F_z$ with the formula $z=z_0 e^{\tau\delta}$. This implies that the flow of the vector field $\delta^{-1}R$ commutes via $h$ with the flow $zd/dz$ on the disc image of $h$. In particular, one can lift the circle with base point $m_0$: $\exp [\sigma]zd/dz\cdot m_0, \sigma\in[0,2i\pi]$. For $\sigma=2i\pi$, we obtain a diffeomorphism $$\rho=\exp [2i\pi]\delta^{-1}R$$ which keep invariant each fiber of $h$. This is the \textit{geometric monodromy} of $\F_h$. The diffeomorphism $\rho$ is periodic with period $\delta$. The orbit of a point $p$ on $F_z$ under the action of $\rho$ is a set of $\delta$ points,
intersection of $F_z$ with the trajectory $T_p$ of $R$ through $p$. The meromorphic first integral of $R$ define a projection $\pi_R$ onto $C$. From the previous description, for any loop $\gamma$ in $F_z$ the $\delta$ elements of its orbit via $\rho$ have the same projection by $\pi_R$ onto a loop which represents an element of $\pi_1'(C,m)$. 
Finally, the elements of $\pi_1'(C,m)$ can be identified to the classes of evanescents loop in a fiber $F_z$ modulo the action of $\rho$, or also to the horizontal family of evanescent loops, obtained by the action of the flow of $R$ on $\gamma$ (the previous description is only the intersection of this family with $F_z$). 

\bigskip

Let $\gamma_1,\cdots \gamma_\mu$ be a basis of the free group $\pi_1(F_z, p)$, and let $\Gamma_1,\cdots \Gamma_\mu$ be their projection in $\pi_1'(C,m)$ or their class modulo $\rho$. We want to compute their image $h_{\Gamma_i}$ in $\hol(\F,\F_h)$ when \F\ has a proper Galois envelope. In the non exceptional case, the probleme reduces to compute $\hol(\F_N,\F_h)$ where $\F_N$ may be defined by the following one-form writen under its adapted form
$$\omega_N=\omega_h-(\sum_{i=1}^{\mu}c_if_i)\delta(h)\frac{\omega_R}{\alpha}.$$
Notice that the $\mu$ one-forms $\eta_i=f_i\frac{\omega_R}{\alpha}$ are horizontal, i.e. invariant under the action of $R$. Indeed, we have
$$L_R(\eta_i)=(i_Rd+di_R)(f_i\frac{\omega_R}{\alpha})=0$$
since $R(f_i)=0$ and $d( \frac{\omega_R}{\alpha})=0$. Let $\eta_c= \sum_{i=1}^{\mu}c_i\eta_i$. The choice of $c=(c_1,\cdots c_\mu)$ completely determine the class of $\F_N$ . Since $\eta_c$ is a horizontal form, its integration $T_i=\int_ {\gamma_i}\eta_c$ only depends on the horizontal classe $\Gamma_i$.

The vector field $\delta(h)R$ is a vertical vector field, i.e. a vector field tangent to the foliation defined by $R$, and its restriction on each fiber doesn't depend on this fiber. If we introduce local coordinates around a point $p$ outside $X$ defined by $s=\int\eta_c$ and $t=h^{1/\delta}$ this vector field is a one variable holomorphic vector field $\theta$ in $t$. In the final final form, we have: $\theta=\frac{t^{q+1}}{1+\lambda t^q}\frac{d}{dt}$.

\begin{thm}\label{hol.rel}
The generators of $\hol(\F_N,\F_h)$ are given by
$$h_{\Gamma_i}(p)=\exp[T_i]\theta\cdot p$$
where the $T_i$'s are the periods of $\eta_c$ on the horizontal cycles $\Gamma_i$.
 \end{thm}

We can remark that this holonomy is an abelian group. This is coherent, since for any solvable group $G$ of germs of diffeomorphisms, the subgroup of its elements which are tangent to the identity is always an abelian one (see \cite{CM} or \cite{P1}). \medskip
 
 \begin{preuve}
 The foliation $\F_N$ is also defined by the vector field $$\alpha X_h+(\sum_{i=1}^{\mu}c_if_i)\delta(h)R$$ or by
  $\frac{\alpha X_h}{f_c}+\delta(h)R$ with $f_c=\sum_{i=1}^{\mu}c_if_i$. The key point here is that the vector field
$\frac{\alpha X_h}{f_c}$ commutes with $\delta(h)R$. Indeed, $[\alpha X_h,R]=0$, $f_c$ is a first integral for $R$ and $\delta(h)$ a first integral for $X_h$. Therefore we have:
$$\exp[\sigma](\frac{ \alpha X_h}{f_c}+\delta(h)R) \cdot p=
\exp[\sigma]\delta(h)R\  \circ \ \exp[\sigma]\frac{ \alpha X_h}{f_c}  \cdot p.$$
If $\sigma$ runs over a segment $[0,T]$ in \C, the first member is a path with origin $p$ into the leaf of $\F_N$ through $p$. Likewise, the term $ \exp[\sigma]\frac{ \alpha X_h}{f_c}  \cdot p$ describe a path of origin $p$ into a leaf of the initial foliation $\F_h$ and $\exp[\sigma]\delta(h)R \cdot q$ is a path into a fiber of the vertical foliation defined by $R$. Therefore, the first member define a lift in $\F_N$ of a path in the initial foliation by the projection defined by $R$. If this path is closed for $\sigma=T$, $\exp[T]\delta(h)R\cdot p$ is its relative holonomy. Since $ \eta_c(\frac { \alpha X_h}{f_c} )=1$, in the (multivalued) coordinate $s=\int \eta_c $  such a flow is a translation, and for the periods $T_i$ of $\eta_c$, the segments $[0,T_i]$ are covering of a basis $\gamma_i$. This proves the theorem.
\end{preuve}

This allows us to characterize the exceptional foliations (i.e. those which have a monogeneous relative holonomy group) on their normal form: 

\begin{cor} $\F_N$ is an exceptional foliation if and only if the quotients of the periods $T_i$ are rational numbers.
\end{cor}

Finally, we can deduce from theorem (\ref{hol.rel}) the following realization theorem:

\begin{thm} We fix the quasi-homogeneous curve $X: h=0$ with Milnor number $\mu$. Given a non exceptional abelian sub-group $H$ of $\diff_1(\C,0)$ generated by $h_1,\cdots h_\mu$, there exits a germ of quasi-homogeneous foliation whose relative holonomy group $\hol(\F,\F_h$) is $H$.
\end{thm}

\begin{preuve} In order to construct the class of \F\, we have to choose  $\alpha X_h+(\sum_{i=1}^{\mu}c_if_i)\delta(h)R$.
Since $H$ is a non exceptional abelian subgroup of $\diff_1(\C,0)$, there exists an analytic vector field $\theta$ and a constants $T_i$ such that $h_i=\exp T_i \theta$ (see \cite{P1}). The vector field $\theta$ induce a unique vector field $\delta(h)R$ whose expression on each fiber of $R$ is $\theta$. We only have to choose the constants $c_i$ which will induce the given relative holonmy. The relationship between the constants $T_i$ and the $c_i$'s is given by
$$T_i=\int_{\Gamma_i}\eta_c=\sum_{j=1}^{\mu}c_j\int_{\Gamma_i} \eta_j.$$
It follows that one should have the matricial equality $T=M\cdot C$ where $T$ is the column of the $T_i$'s,
$C$ is the column of the $ c_i$'s and $M=(m_{i,j})$ with $m_{i,j}=\int_{\Gamma_i}\eta_j$. The coefficients of this matrix are constants since $\eta_j$ and $\Gamma_i$ are horizontal. Since the loops $\gamma_i$ generate a basis of the homology of the Milnor fiber and the $\eta_j$ a basis of its cohomology, it is a well known fact that this matrix is an inversible one. Therefore, we may compute $C$ from $T$.
\end{preuve}

\textbf{Remark.} If we admit that the relative holonomy group is a complete invariant of the foliation, the previous result gives us another proof of the main implication $(4)\Rightarrow (3)$ in (\ref{Red.gal.nonex}). Indeed, if \F\ is a Liouvillian foliation, its relative holonomy group is a non exceptional abelian subgroup $H$ of $\diff_1(\C,0)$. We can realize it by another foliation given by a normalized vector field $X_N$, whose algoritmic invariant $\mathcal {L}(\F_N)$ is one dimensional. Since the two foliations are analytically equivalents, we have $\mathcal {L}(\F_N)=\mathcal {L}(\F)$ up to a conjugacy, and $\mathcal {L}(\F)$ is also one-dimensional.

\bigskip
\section{Open problems}

\medskip

In the present class of quasi-homogeneous foliations, there remain the following questions:
\begin{itemize} 
\item Find the relation between the algorithmic invariant $\mathcal {L}(\F)$ and the geometric one (relative holonomy) outside the Galois reducible case. In the general case, this transcendental relation will not reduce to the exponential of one variable vector field. Probably, we shall have to consider Campbell-Hausdorff type formulae;
\item Prove the $k$-summabuility for the final normal forms and find the geometric meaning of this order $k$.
\end{itemize}
One can try to extend such a study to any germ of foliation in $\C^2$: 
\begin{itemize} 
\item In the non dicritical case (i.e. when the exceptional divisor is an invariant set of the foliation), outside the quasi-homogeneous context, we have no normal forms. We would like to construct them, having in mind the present motivations: a good representative of an holomorphic foliation may allow us to compute its Galois closure, and its geometric invariants. Of course, we agree divergent models in order to get the previous conditions, and we expect their summability.
\item In the generic dicritical case (i.e. when the foliation is desingularized after one blowing up such that the projective line is not an invariant set), we have formal normal forms: see \cite{ORV}. Can we make use of these models to compute their Galois closure?
\end{itemize}
We can also consider the following developments:
\begin{itemize} 
\item (suggested by B. Malgrange) study the Galois closure for any local codimension one foliations: can we reduce it to the previous dimension two cases? 
\item study the Galois closure of vector fields in $(\C^2,0)$. This means that we first have to classify vector fields not only up to a unity, and to construct formal normal forms with respect to this classification. 
\item develop a similar study for an algebraic foliation on the projective plane near an algebraic invariant set.
\end{itemize}

\begin{center}
-----------------------------
\end{center}

\begin{thebibliography}{99}

\bibitem{CDS} M. Canalis-Durand, R. Schaefke -
\textit{Divergence and summability of normal forms of systems of differential equations with nilpotent linear part.}  
Ann. Fac. Sci. Toulouse Math. (6)  13  (2004),  no. 4, 493--513. 

\bibitem{CAS0} G. Casale -
\textit{Suites de Godbillon-Vey et intégrales premières}
C.R. Acad. Sci. Paris 335 (2002). 

\bibitem{CAS1} G. Casale - 
\textit{$D$-enveloppe d'un difféomorphisme de $({\Bbb C},0)$. }
 Ann. Fac. Sci. Toulouse Math. (6)  13  (2004),  no. 4, 515--538.

\bibitem{CAS2} G. Casale - 
\textit{D-enveloppe d'une application rationnelle,}
Publicaciones mathemàtiques (Barcelona) 50-1 (2006)

\bibitem{CAS3} G. Casale - 
\textit{Feuilletages singuliers de codimension un, groupoïde de Galois et intégrales premières,}
Ann.  Inst. Fourier 56-3 (2006)

\bibitem{CAS4} G. Casale - 
\textit{Le groupoïde de Galois de l'équation de Painlevé I et son irréductibilité,} to appear in Comm. Math. Helv.

\bibitem{CM} D. Cerveau and R. Moussu -
\textit{Groupes d'automorphismes de $(\C,0)$ et équations diffé\-ren\-tielles $y\ dy+\cdots =0.$ }
Bull. Soc. math. France {\bf 116} (1988), 459--488.

\bibitem{YG} Y. Genzmer -
\textit{Phd thesis,} lab. E. Picard Toulouse (2006)

\bibitem{LOR} F. Loray -
\textit{Feuilletages holomorphes à holonomie résoluble}
Phd Thesis, University of  Rennes I (1994)

\bibitem{LOR1} F. Loray - 
\textit{Réduction formelle des singularités cuspidales de champs de vecteurs analytiques.}
J. of Diff. Equations {\bf 158, 1} (1999), 152--173.

\bibitem{LM} F. Loray and R. Meziani - 
\textit{Classification de certains feuilletages associés à un cusp.}
Bol. Soc. Bras. Math. {\bf 25} (1994), 93--106.

\bibitem{MAL0} B. Malgrange -
Notes manuscriptes (2000).

\bibitem{MAL1} B. Malgrange -
\textit{Le groupoïde de Galois d'un feuilletage,} 
Monographie \textbf{38}, vol. \textbf{2} de l'Enseignement Mathématiques (2001).

\bibitem{MAL2} B. Malgrange -
\textit{On the non linear Galois theory,}
Chinese Ann. Math. ser B, \textbf{23, 2}, (2002)

\bibitem{MR1} J. Martinet and J.P. Ramis -
\textit{Classification analytique des équations diffé\-ren\-tielles non linéaires résonnantes du pemier ordre.}
Ann. Sci. Ecole Norm. Sup., t.{\bf 16} (1983), 571--621.

\bibitem{JFM} J.F. Mattei -
\textit{Quasihomogénéité et équiréductibilité de feuilletages holomorphes en dimension 2,}
Astérisque {\bf 261} S.M.F. (2000) 253--276.

\bibitem{MM} J.F. Mattei and R. Moussu -
\textit{Holonomie et intégrales premières,}
Ann. Sci. Ecole Norm. Sup., t.{\bf 13} (1980), 469--523.

\bibitem{M} R. Moussu -
\textit{Holonomie évanescente des équations différentielles dégénérées transverses.}
Singularities and Dynamical Systems, S.N. Pnevmatikos ed., Elsevier Sc. Publisher (1985) 161--173

\bibitem{ORV} L. Ortiz-Bobadilla, E. Rosales-Gonz\' alez, V.M. Voronin -
\textit{Rigidity theorems for generic holomorphic germs of dicritic foliations and vector fiels in $(\C^2,0)$,}
Mosc. Math. J 5 (2005), no 1, 171--206.

\bibitem{P1} E. Paul -
\textit{Feuilletages holomorphes singuliers \`a holonomie r\'esoluble, }
J. reine angew. Math. {\bf 514} (1999), 9-70.

\bibitem{P2} E. Paul -
\textit{Formal normal form for a perturbation of a quasi-homogeneous hamiltonian vector field.}
Journal of Dynamical and Control Systems {\bf 10} (2004) 545-575.

\bibitem{SEI} A. Seidenberg -
\textit{Reduction of singularities of the differential equation $Ady=Bdx$,}
Amer. J. Math. \textbf{90} (1968)

\bibitem{ZOL} E. Str\'ozyna, H. Zoladek - 
\textit{The analytic normal form for the nilpotent singularity.}
J. of Diff. Equations {\bf 179, 2} (2002), 479--537.




\end{thebibliography}
\end{document}